\documentclass[12pt]
{article}
\usepackage{amssymb,amsmath}

\usepackage[latin1]{inputenc}
\usepackage[T1]{fontenc}
\usepackage{play}

\usepackage{graphicx}
\DeclareGraphicsExtensions{.pdf,.png,.jpg}

\graphicspath{ {images/} }

\usepackage{lipsum}

\let\OLDthebibliography\thebibliography
\renewcommand\thebibliography[1]{
  \OLDthebibliography{#1}
  \setlength{\parskip}{0pt}
  \setlength{\itemsep}{0pt plus 0.3ex}
}

\usepackage{amsmath}
\usepackage{amsthm}
\usepackage{amsfonts}
\usepackage{bbm}
\usepackage{amssymb}
\usepackage{graphicx}
\usepackage{color}

\usepackage{wasysym, marvosym}

\usepackage{fancybox}
\newlength{\querylen}
\newcommand{\prob}{\mathbb{P}}

\newcommand{\E}{\mathbb{E}}

\theoremstyle{definition}

\theoremstyle{remark}

\begin{document}
\title{On sequential selection  and a first passage problem for the Poisson process} 
\author{Alexander Gnedin\thanks{Queen Mary, University of London}}
\maketitle
\begin{abstract}
\noindent
This note is motivated by  connections between the online and offline problems of selecting a possibly long subsequence from a Poisson-paced  
sequence of uniform marks under either a monotonicity or a sum constraint. 
The offline problem with the sum constraint  amounts to counting the Poisson arrivals before their total exceeds a certain level.
A precise asymptotics  for the mean count  is obtained by  coupling with  a  nonlinear pure birth process.

\end{abstract}

\noindent
\section{Introduction}

{\it When a shuttle  carrying a large number of hotel guests arrives at the hotel,
 the passengers start queuing and  pass the exit door at times of a Poisson process.  
The waiting times  spent in the queue are added up as the passengers quit.
What is the number $N(t)$ of passengers that  exit the shuttle before the accumulated waiting time exceeds $t$?}

We shall call this the shuttle exit problem.  The {\it exit count} $N(t)$ is the maximum number of Poisson times with the total not exceeding $t$. 
The total waiting time and the 
 exit count process are important in many models of applied probability. 
Our interest stems from the
 connection to the online version of  the longest increasing subsequence problem with Poisson arrivals, which we now describe.

Suppose  independent, uniform [0,1] marks  arrive sequentially  at times of a unit rate  Poisson process on $[0,t]$. A prophet 
with complete overview of the data can use an offline algorithm to
 select the longest increasing 
subsequence of  length $L^*(t)$.   
A nonclairvoyant   gambler learns  the data  and makes irrevocable decisions in real time using a nonanticipating online selection policy. Let $L(t)$ be the length of increasing subsequence selected under the 
online policy that achieves
the maximum expected length.
As $t\to\infty$,
\begin{eqnarray}\label{LIS} 
\E L^*(t) &= & 2\sqrt{t} -c_0\,t^{1/6}+o(t^{1/6}) ,~~\\
\label{LOIS}
\E L(t)~ &=  & \sqrt{2t}-\tfrac{1}{12}\log t +O(1),
\end{eqnarray}
(where  $c_0=1.77\ldots$ is an explicit constant).
The limit ratio  $2:\sqrt{2}$  serves as  a rough measure of  advantage of the prophet over the gambler. 
The  asymptotics (\ref{LIS})   has a long and colourful history, culminating in the work by 
 Baik, Deift and Johansson \cite{BDJ}. See Romik's book \cite{Romik} for a nice exposition.
The leading term 
of (\ref{LOIS}) is due to Samuels and Steele \cite{SS} 
who were first to study the online problem, later Bruss and Delbaen \cite{BD} identified the logarithmic order of the second term
and the full expansion has appeared  recently in \cite{GSBernoulli}.

Remarkably, the online increasing subsequence problem can be recast as a  very different stochastic task, with
the monotonicity constraint  replaced by the condition that the sum of selected marks should not exceed $1$. The latter is commonly interpreted as a bin-packing problem, where gambler's objective is  to maximise the expected number of items packed online in a bin of unit capacity \cite{BK, Coffman}. 
By analogy with (\ref{LIS}) and (\ref{LOIS}) it is natural to consider the offline counterpart of $L(t)$
in the bin packing context.
Obviously, with full information,
the optimal prophet's policy amounts to the {\it smallest first policy} that 
packs the items in the increasing order of size as long as they fit in the bin.

Since the marks sorted into increasing order themselve comprise a homogeneous Poisson process,
zooming in  the marks scale with factor $t$ and  changing the metaphore, 
it is seen that the number of  items packed under the smallest first policy coincides with the exit count $N(t)$ from the shuttle problem we started with.

The first surprise in the online-offline bin packing comparison  comes with the fact that the limit prophet-to-gambler ratio is equal to $1$.
This follows from the asymptotics
$\E N(t)\sim \sqrt{2t},$
which in turn can be concluded from a benchmark   \cite{BK, BR, Gnedin, SteeleBR} upper bound $\E N(t)<\sqrt{2t}$, the trivial inequality 
$L(t)\leq N(t)$ 
 and  (\ref{LOIS}). 
Therefore to assess the magnitude of  prophet's advantage one needs to examine the finer the mean exit count more closely.

In this note we find a formula for $\E N(t)$ in terms of the Borel distribution. Though explicit, the formula  seem to  require substantial analytic work to extract  the desired second term of the asymptotic expansion.
We circumvent this
 by resorting to elementary probabilistic tools, with 
the core of our approach being the observation that $N(t)$, for each fixed $t$, has the same distribution as the {\it entrance count} $M(t)$ appearing in the following dual 
shuttle entrance problem.

{\it When the shuttle picks up hotel guests at the airport, they enter by the Poisson process. The shuttle departs at the moment when  the total waiting time of the driver and all passengers inside 
the shuttle  is $t$. What is  the number $M(t)$  of hotel guests in the shuttle  by the departure?}

We observe that the process $M(t)$ is a nonlinear  pure-birth Markov chain which was considered in Kingman and Volkov  \cite{Kingman} in the context
of gunfight models. 
Using the identity in distribution we show that
\begin{equation}\label{twothirds}
\sqrt{2t}-\E N(t)\to \frac{2}{3}
\end{equation}
and that the difference is always less than 1.
This contrasts sharply with the second terms in (\ref{LIS}) and (\ref{LOIS}).
For   
the difference 
between the prophet and gambler values  we have therefore
$$\E N(t)-\E L(t)\sim \tfrac{1}{12}\log t.$$

Bruss and Delbaen \cite{BD} showed that $L(t)$ is ${\rm AN}(\sqrt{2t},\tfrac{1}{3}\sqrt{2t})$ (AN abbreviates `asymptotically normal'), see also \cite{GSBernoulli}.
We argue that the same is true for $N(t)$. The asymptotic  coincidence of variances looks unexpected since the underlying selection policies are very different.
We remind that in  the increasing subsequence problem the  {\it types} of the limit  distribution of $L^*(t)$ and $L(t)$ are different, as
the distribution of the maximum offline length $L^*(t)$ approaches the Tracy-Widom law  from the random matrix theory \cite{BDJ, Romik}.

This note is a collection of  snapshots around (\ref{twothirds}).
To keep the discussion short, details of routine proofs are only sketched.
Related work on  sums of  consequitive arrivals  in the case of inhomogeneous rate appeared in \cite{AKZ},
and on the integrated Poisson process  in  \cite{Suyono}.

The rest of the paper is organised as follows.
In the next two sections we add insight to what is already known regarding the coupling of online selection problems and the benchmark upper bound. 
In sections 4 and 5 we scrutinise the exit-entrance duality.
In section 6 we record the normal limits.
In section 7 we derive a series formula for the mean count. 
In section 8 we employ the pure birth process to refine the $\sqrt{2t}$ asymptotics.
In section 9  we depoissonise (\ref{twothirds}) to improve upon the well known \cite{BK, BR, Coffman, SteeleBR} fixed sample asymptotics of the smallest first policy. 
A large deviation bound needed for our arguments is derived in the last section.

Throughout we shall be using the notation 
$$\nu(t):=\E N(t), ~~\sigma^2(t):={\rm Var} N(t).$$

\section{Coupling of online problems}

We first detail the equivalence between the online increasing subsequence and bin packing problems.
The question about  explicit coupling was  emphasized in Section 5 of Steele \cite{SteeleBR}, where problems with fixed number of arrivals $n$ were discussed.

The distribution of marks in the increasing subsequence problem does not matter (subject to being continuous), 
while the bin packing problem is not distribution-free.
In the special case of uniform $[0,1]$ marks and  the bin of unit capacity,
the equivalence  in terms of the optimal policies 
is commonly argued by comparing the dynamic programming equations for the value function \cite{Arlotto, Coffman}.
It is also noticed in  \cite{Coffman} (p. 455) that the greedy online bin packing policy translates as the increasing sequence of record marks.

The following construction provides  a general coupling  in our setting with the Poisson arrivals, but 
it can be readily adjusted to other arrival processes including the discrete time models with fixed
or random horizon  
\cite{Arlotto, Gnedin, SS}.

Let $\Pi^{\tiny\Square}$ be a planar Poisson point process with unit rate in  the strip $[0,\infty)\times[0,1]$. We  endow $\Pi^{\tiny\Square}$   with the natural filtration generated by 
$\{{\Pi^{\tiny\Square}}|_{[0,t]\times[0,1]}, ~t\geq 0\}$.
The generic atom  of $\Pi^{\tiny\Square}$  at location  $(\tau,\xi)$ is understood as mark $\xi$ arriving at time $\tau$.

We define an {\it  $i$-selection policy} to  be a nondecreasing, adapted,  c{\'a}dl{\'a}g  jump  process   $I$ with $I(0)=0$, such that the north-west corners of the graph of $I$ 
are some atoms $(\tau_k,\xi_k)$ of $\Pi^{\tiny\Square}$ labeled 
by increase of the time component.
This sequence of atoms spanning the graph is an increasing chain in the partial order in two dimensions.

Similarly, we define a {\it  $b$-selection policy}  to  be a nondecreasing, adapted,  c{\'a}dl{\'a}g  jump  process   $B$ with $B(0)=1$ and values in $[0,1]$.
We require  that each jump be corresponding to  an atom  $(\tau_k,\xi_k)$,  so that the jump-time is $\tau_k$ and the increment is $\xi_k$.
Thus the range of $B$ is the sequence of partial sums of $\xi_1, \xi_2,\dots$.

For a fixed $i$-selection policy $I$, we are going to introduce an invertible  random transform 
$\phi_I$ of $[0,\infty)\times[0,1]$, which will map $I$ to  a $b$-selection policy with the same path $B=I$.
The construction is iterative.

At each step $k$ we shall have  $[0,\infty)\times[0,1]$
and its  duplicate obtained  by a measure-preserving $\beta_k$.
Start with two identical copies of the strip
equipped with Poisson point scatters of $\Pi^{\tiny\Square}$,  and a fixed path of $I$ spanned on some points $(\tau_k,\xi_k)$.
Let $\beta_0$ be the identity, and $\xi_0=0$.
At step $k>0$ only the strip $\beta_k((\tau_k,\infty)\times [\xi_{k-1},1])$ undergoes a change  which amounts to  cutting   at height $\xi_k-\xi_{k-1}$  by the horizontal line and 
placing 
part $\beta_k((\tau_k,\infty)\times [\xi_{k-1},\xi_k])$  atop of $\beta_k((\tau_k,\infty)\times [\xi_k,1])$  with the orientation preserved.
The mapping $\beta_{k+1}$ is the composition of $\beta_k$ and this surgery.  
With probability one, each point moves under $\beta_k$'s finitely many times, as the moves may only be associated with $(\tau_k,\xi_k)$'s to the left of this point. 
Thus we may define $\phi_I$ as the composition of all $\beta_k$'s.

Note that $\phi_I$  preserves the planar Lebesgue measure and 
does not alter the time component, so leaving each set $(t,\infty)\times[0,1]$ invariant.
Consider the transformed point process $\widehat{\Pi}^{\tiny\Square}:=\phi_I({\Pi}^{\tiny\Square})$.
By the invariance, ${\Pi}^{\tiny\Square}$  and 
$\widehat{\Pi}^{\tiny\Square}$ share the same one-dimensional Poisson process of arrival times. 
Given arrival at time $\tau$, the image of $(\tau,\xi)$ under $\phi_I$ is uniformly distributed on $\{\tau\}\times[0,1]$ and is independent of ${\Pi}^{\tiny\Square}|_{[0,\tau)\times[0,1]}$, hence also independent of
$\widehat{\Pi}^{\tiny\Square}|_{[0,\tau)\times[0,1]}$. But this implies that  $\widehat{\Pi}^{\tiny\Square}$ has the same distribution as ${\Pi}^{\tiny\Square}$.
The transformation $\phi_I$ sends the sequence $(\tau_k,\xi_k)$  to  a sequence $(\tau, \xi_k-\xi_{k-1})$ (where $\xi_0=0$),  which are now some
 atoms of $\widehat{\Pi}^{\tiny\Square}$,  and $I$ becomes a
$b$-selection policy spanned on the transformed sequence.

The above concepts of selection policy are much  more general than the Markovian threshold policies    studied in the literature. For the purpose of optimisation, however, 
it is sufficient to consider  the following family  of policies.
For  $\psi:[0,\infty)\to [0,1]$ thought of as a function controlling the size of acceptance window, 
and given horizon $t$, an $i$-selection policy is defined recursively by the  rule:
conditionally on 
arrival occurring at time $\tau<t$
and  given $I(\tau-)=x$ (the last selection so far),
the observed mark $\xi$ is selected if and only if 
\begin{equation}\label{AC1}
0<\frac{\xi-x}{1-x}\leq\psi((t-\tau)(1-x)).
\end{equation}
In \cite{GSBernoulli} we called such policies self-similar because the performance from each stage on
only depends on the mean number of future acceptable arrivals. 
Thus defined, $I$ is a jump Markov process with transition mechanism determined by $\psi$.
The twin $b$-selection policy has the acceptance condition 
\begin{equation}\label{AC2}
0<\frac{\xi}{1-x}\leq\psi((t-\tau)(1-x)),
\end{equation}
given $B(\tau-)=x$ (the total of selected items so far).  
The optimal  $i$-/$b$-selection policy is of this form with some control $\psi^*$ satisfying 
$$\psi^*(z)\sim \sqrt{\frac{2}{ z}}- \frac{1}{3z}, ~~~~z\to\infty,$$ 
 see \cite{BG, BD, GSBernoulli}.
In \cite{GSBernoulli} we proved that every policy having $\psi(z)\sim \sqrt{2/z}$ 
is within $O(1)$ from the optimality, that is achieves the asymptotics (\ref{LOIS}).

The general Markovian policy differs from (\ref{AC1}) and (\ref{AC2}) in that $\psi$ is replaced by the general function of $\tau, t$ and $x$.
Notable other examples are the {\it greedy} policy with the function $1$ and the {\it stationary} policy with the function $\sqrt{2/t}\wedge 1$.

\section{The upper bound}

For the rest of this paper the variable $t$ will have the meaning of  either the bin capacity (the offline bin-packing contest) or the total waiting time (the shuttle context).
For the time parameter of the Poisson process we shall use the variable $x$.

Let  $\pi_1<\pi_2<\dots$ be the points of a unit rate Poisson process $\Pi$ on the positive half-line. 
The exit count is defined as     
$$N(t):=\max\{n: \pi_1+\dots+\pi_n\leq t\},~~~~t\geq 0,$$
 where $\max\varnothing=0$.

There is a benchmark  upper bound  for the mean,
\begin{equation}\label{sqrt}
\nu(t)<\sqrt{2t}\,, ~~~~~~~t>0,
\end{equation}
that appeared in the Poisson setting in \cite{BK} (Example 2.4).
Similar inequalities for sums of order statistics from the general distribution are found in
\cite{BR}, also 
see \cite{SteeleBR} for extended discussion. 
We relate (\ref{sqrt}) to an isoperimetric inequality, much in line with the examples from \cite{BG, Gnedin}.

Fix $t$. The set of Poisson points $\pi_n$ with $\pi_1+\dots+\pi_n\leq t$ is a point subprocess
of $\Pi$ with rate function $p_t$ 
 satisfying
\begin{equation}\label{Iso}
\nu(t)=\int_0^t \, p_t(x) {\rm d}x,~~~ \int_0^t x\, p_t(x) {\rm d}x\leq t, ~~~~0\leq p_t(x)\leq 1.
\end{equation}
This suggests a problem from  the calculus of variations,
$$
\int_0^t \, q(x) {\rm d}x\to\max ,~~~ \int_0^t x\, q(x) {\rm d}x\leq t, ~~~~0\leq q(x)\leq 1.
$$
The Lagrangian function becomes
$$\int_0^t (\theta-x)q(x){\rm d}x, ~~~{\rm with~~} \theta>0,$$
which for given multiplier $\theta$ is maximised by the indicator function $q(x)=1(x\leq \theta)$.
Accounting for the constraint, the overall maximum value of the integral is $\sqrt{2t}$, 
 attained at 
$$\theta^*=\sqrt{2t},~~~  q^*(x)=1(x\leq\sqrt{2t}),$$
which gives the upper bound
 (\ref{sqrt}) follows.

\noindent
{\bf Remark} Solution $q^*$ corresponds to a packing policy  that picks all items smaller than the threshold $\sqrt{2/t}$.
The policy 
violates the (almost sure) sum constraint but meets 
a weaker mean-value constraint. This policy is online implementable and outputs the number of selections with Poisson$(\sqrt{2t})$ distribution,
so has the variance about three times higher than under the optimal offline (see below) or the optimal online policy \cite{BD, GSBernoulli}.

\section{The exit-entrance duality}

Consider the shuttle entrance problem.
When the $n$th passenger enters 
the total waiting time of everyone inside the shuttle is
$$\pi_1+2(\pi_2-\pi_1)+\dots+n(\pi_n-\pi_{n-1})=n\pi_n-(\pi_1+\dots+\pi_{n-1}),$$
so the entrance count is 
$$M(t):=\max\{n: n\pi_n-(\pi_1+\dots+\pi_{n-1})\leq t\}.$$

We assert that 
\begin{equation}\label{N=M}
N(t)\stackrel{d}{=}M(t).
\end{equation}
Indeed, since 
$$N(t)\geq n \Leftrightarrow \pi_1+\dots+\pi_n\leq t, ~~~M(t)\geq n \Leftrightarrow n\pi_n- (\pi_1+\dots+\pi_{n-1})\leq t,$$
we need to check that 
$$
\pi_1+\dots+\pi_n\stackrel{d}{=}n\pi_n- (\pi_1+\dots+\pi_{n-1}).
$$
Recall that, given $\pi_n$,  the quotients $\pi_j/\pi_n, ~j<n,$ are independent from $\pi_n$ and jointly distributed like the uniform order statistics. Thus
for $u_1,u_2,\dots$ iid uniform $[0,1]$  we have
\begin{eqnarray*}
\pi_1+\dots+\pi_n=\pi_n\left(1+ \frac{\pi_1+\dots+\pi_{n-1}}{\pi_n} \right)
\stackrel{d}{=}\pi_n(1+u_1+\dots+u_{n-1}) \stackrel{d}{=}\\
\pi_n(1+(1-u_1)+\dots+(1-u_{n-1})) \stackrel{d}{=}
\pi_n(n-(u_1+\dots+u_{n-1}))\stackrel{d}{=}\\
\pi_n\left(n-\frac{\pi_1+\dots\pi_{n-1}}{\pi_n}\right)=n\pi_n- (\pi_1+\dots+\pi_{n-1}),
\end{eqnarray*}
where we used symmetry of the uniform distribution.

It is also instructive to argue in terms of the iid exponentially distributed gaps $\eta_j:=\pi_j-\pi_{j-1}$ (with the convention $\pi_0=0$). 
We have 
\begin{eqnarray*}
\pi_1+\dots+\pi_n&=&n\eta_1+(n-1)\eta_2+\dots+\eta_n\stackrel{d}{=}\\
 \eta_1+2\eta_2+\dots+n\eta_n&=&\pi_1+2(\pi_2-\pi_1)+\dots+n(\pi_n-\pi_{n-1}).
\end{eqnarray*}
The variables $\zeta_n:=\eta_1+2\eta_2+\dots+n\eta_n$ are the jump-times of the entrance count process.  Thus $(M(t),~t\geq 0)$ is a pure-birth process that starts with
$M(0)=1$ and moves from state $n$ to state $n+1$ at rate $(n+1)^{-1}$.

The entrance count process has a simple combinatorial  interpretation. 
Think of an urn with one red  and some number of white balls. At times of the Poisson process a ball is randomly chosen and replaced to  the urn.
If the chosen ball is red, a white ball is added to the urn, otherwise the urn composition is not changed.
For the process starting with one red ball, $M(t)$ is the number of  white balls in the urn at time $t$.

The identity  (\ref{N=M}) only holds for the marginal distributions, and the exit count  process $(N(t),~t\geq0)$ is not even Markovian.
The driver's waiting time was included in the total waiting time to avoid a shift in the distributional identity.
We note in passing that without appealing to   (\ref{N=M}) the upper bound
$\E M(t)\leq \sqrt{2t}$
does not seem at all obvious.

\section{Integrals of the Poisson process}

Let 
$$
T(x) :=   \int_0^x y\, {\rm d}\Pi(y)= \sum_{j=1}^{\Pi(x)}\pi_j   , ~~~S(x):= \int_0^x\Pi(y){\rm d}y= \sum_{j=1}^{\Pi(x)} (x-\pi_j) .
$$
The total waiting time accumulated within the real time $x$ is $T(x)$ in the shuttle exit problem, and $x+S(x)$ in the entrance problem, where
$x$ is added to account for driver's waiting time.
The integration by parts formula becomes
$$T(x)=x\Pi(x)-S(x).$$

By reversibility of $\Pi$ on $[0,x]$ we have
\begin{equation}\label{T=S}
T(x)\stackrel{d}{=}S(x),
\end{equation}
This identity has appeared  in \cite{Suyono}, where it was concluded analytically from the identity of   Laplace transforms. 
Despite that (\ref{T=S}) holds for each fixed $x$, the processes are very different: $T$ is a jump process with independent increments, while the paths of $S$ are  piecewise linear.

Plugging for $x$  the Poisson times we obtain a few  `total waiting time paradoxes'. 
First note the defining recursions 
\begin{equation}\label{defrec}
S(\pi_{n+1})=S(\pi_n)+n(\pi_{n+1}-\pi_n),~~ T(\pi_{n+1})=T(\pi_n)+\pi_{n+1}.
\end{equation}
Now,
given $\pi_{n+1}$, the variables $T(\pi_n)$ and $S(\pi_{n+1})$
have the same distribution, and so  unconditionally  
\begin{equation}\label{TpiS}
T(\pi_n)\stackrel{d}{=} S(\pi_{n+1}),
\end{equation}
in apparent disagreement with  (\ref{T=S}).
Moreover,
$S(\pi_{n+1})\stackrel{d}{=} S(\pi_n)+\pi_n$, which is to be compared with  (\ref{defrec}) and (\ref{TpiS}).
The latter identity is equivalent to 
$$\pi_{n+1}(u_1+\dots+u_n)\stackrel{d}{=}\pi_n(1+u_1+\dots+u_{n-1}),$$
where the $\pi_n$'s are independent of the iid uniform $u_j$'s. 
To prove the last formula directly, one can observe two ways to   split $T_n$ in independent factors,  as
 $\pi_{n+1}(T(\pi_n)/\pi_{n+1})$ and $\pi_n(T(\pi_n)/\pi_n)$, then
represent the quotients in brackets
 in terms of the $u_j$'s. 
See \cite{GM} for more involved exponential-uniform  identities derived from the  planar Poisson process.

Next, we aim to represent the exit and entrance counts as time-changed Poisson process.
Let $X(t)=\min\{x: T(x)>t\}$ be the right-continuous inverse of $T$, with $X(0)=\pi_1$. 
We can take here $\min$ rather than infinum since $T$ jumps at the discrete set of Poisson points.
We have then
\begin{equation}\label{NPi}
N(t)=\Pi(X(t))-1.
\end{equation}
The process  $S(x)+x$ is strictly increasing, so there is a well defined inverse $\tau$ with $S(\tau(t))+\tau(t)=t$ and
$$\frac{{\rm d} \tau}{{\rm d} t}=\frac{1}{\Pi(\tau(t))+1}.$$
The entrance counting process satisfies
\begin{equation}\label{MPi}
M(t)=\Pi(\tau(t)).
\end{equation}
The last two formulas give yet another proof that the entrance count is a pure-birth process with the jump rate $(n+1)^{-1}$ at state $M(t)=n$.

\section{Normal limits}

Note that $T$ has independent increments.
Application of  Campbell's formula  yields the
 moments 
$$\E\, T(x)=\tfrac{1}{2}x^2, ~~~{\rm Var}\, T(x)=\tfrac{1}{3}x^3,$$
and, more generally, the moment generating function
$$\E e^{zT(x)}= \exp\left( \frac{e^{zx}-zx-1}{z}\right).$$
Inverting this, Suyono and  van der Weide \cite{Suyono} found the density of $T(x)$ in terms of modified Bessel functions
(note that $T(x)$ has mass $e^{-x}$ at zero).

Routine application of the law of large numbers and the central limit theorem show that for $x\to\infty$ 
\begin{eqnarray*}
T(x), S(x)\sim \tfrac{1}{2}x^2 {\rm ~~a.s.}, ~{\rm and~are}~~{\rm AN}\left( \tfrac{1}{2}x^2,\tfrac{1}{3}x^{3}\right).
\end{eqnarray*}
Inverting these asymptotic relations in a way familiar from the renewal theory,  using (\ref{NPi}), (\ref{MPi}) and the asymptotics of $\Pi$ itself,
 we obtain for $t\to\infty$ that
\begin{eqnarray}\label{CLT-NM} 
N(t), M(t) \sim \sqrt{2t}~{\rm a.s.}, ~{\rm and~are}~~{\rm AN}\left( \sqrt{2t}\,, \tfrac{1}{3}\sqrt{2t}\right).
\end{eqnarray}

The representation 
$$M(t)= \max\{n: \zeta_n\leq t\}$$
embeds the analysis of the `renewal function' $\nu(t)=\E M(t)$ into the general framework of the
renewal theory with nonhomogeneous inter-arrival times \cite{Smith}.
Let 
$$a_n: = \tfrac{1}{2}\,{n(n+1)}, ~~b_n: =\tfrac{1}{6}\,{n(n+1)(2n+1)}.$$
It is routine to see that $\E\, \zeta_n = a_n, ~{\rm Var}\, \zeta_n =b_n^2$
and that, 
$$\zeta_n\sim a_n~~{\rm a.s.}, ~~~{\rm ~ and~is~~}
{\rm AN}(a_n, b_n^2).$$
Inverting this  yields another,  more straightforward,   proof of (\ref{CLT-NM}).

The normal limit suggests the asymptotics for the variance
\begin{equation}\label{asympSigma}
\sigma^2(t)\sim \tfrac{1}{3} \sqrt{2t}.
\end{equation}
For a time being we shall take the formula for granted, deferring its 
justification, by checking the uniform integrability, to the last section of this paper.

\section{Exact formulas}

Recall that $\zeta_n=\sum_{j=1}^n j\eta_j$ (with the $\eta_j$'s being  iid exponential), which has the same distribution as the entrance total waiting time $S(\pi_n)+\pi_n$.

The Laplace transform of $\zeta_n$ is
$${\mathbb E}\exp(z\zeta_n)=\prod_{j=1}^n \frac{1}{1-jz}.$$
Inverting this 
yields a formula for the distribution function 
$${\mathbb P}(M(t)\geq n)={\mathbb P}(\zeta_n\leq t)= \frac{1}{n!}\sum_{j=1}^n {n\choose j}(-1)^{n-j} j^n (1-e^{-t/j}).$$
See \cite{Van} for an asymptotic expansion for large $t$.

For the mean of $N(t)$, with a small series work, we obtain an exact formula
\begin{equation}\label{exact}
\nu(t)=\sum_{n=1}^\infty {\mathbb P}(\zeta_n\leq t)=\sum_{j=1}^\infty e^{-j}\frac{j^j}{j!}(1-e^{-t/j}).
\end{equation}

Intriguingly,  (\ref{exact}) can be viewed as  a mean  over the Borel distribution
$${\mathbb P}(Z=j)=e^{-j} \frac{j^{j-1}}{j!},$$
which is the law for the total offspring in the branching process with the Poisson$(1)$ reproduction. Specifically,
\begin{equation}\label{Borel}
\nu(t)={\mathbb E}(Z(1-e^{-t/Z}))= t \int_0^\infty  e^{-ty}\, {\mathbb E}\,(X \,\, 1(X\leq y^{-1})){\rm d}y.
\end{equation}

\section{Bounds on the mean and the limit constant}

The  transition probability of the entrance count process $M$ is
$$\prob(M(t+dt)-M(t)=1|M(t))=\frac{dt}{M(t)+1},$$
which upon taking the expectation becomes
\footnote{Consider the general pure-birth Markov chain $M$ with $M(0)=0,$  and transition rates 
$\beta_n,\,n\geq 0,$ meeting the regularity condition 
$\sum_{n=0}^\infty \tfrac{1}{\beta_n}=\infty$.
The mean population size $\nu(t):=\E M(t)$ and the second moment satisfy then $\nu'(t)=\E \beta_{M(t)}$ and
$(\E M^2(t))'=\E \,[(2M(t)+1)\beta_{M(t)}]$.
}
\begin{equation}\label{eqnnu}
\nu'(t)= \E \left(\frac{1}{M(t)+1}\right).
\end{equation}
Applying Jensen's inequality we arrive at a differential inequality
$$\nu'(t) > \frac{1}{\nu(t)+1}\,,$$
which is readily solved by separating variables as
$\nu(t)>\sqrt{2t+1}-1,~~~t>0.$

So together with (\ref{sqrt}) we  have fairly tight bounds
\begin{equation}\label{bounds}
\sqrt{2t+1}-1<\nu(t)<\sqrt{2t}, ~~~~~t>0,
\end{equation}
where the gap
stays below  1 for all $t$. 
The bounds (\ref{bounds}) clearly suggest that the gap converges to a constant.

Next, we  aim at finding the constant perceived from  (\ref{bounds}).
The random variable 
$$M(t)+1=\min\{n: \zeta_n>t\}$$
is a stopping time.
Doob's optional sampling theorem applied  to the martingale $\zeta_n-a_n$ yields a Wald-type identity
\begin{equation}\label{Wald}
\E \left[\zeta_{M(t)+1}-\tfrac{1}{2}(M(t)+1)(M(t)+2)\right]=0.
\end{equation}
On the other hand, conditionally on $M(t)=n$ the distribution of $\zeta_{M(t)+1}-t$ is exponential with rate $(n+1)^{-1}$, so unconditionally
we can write the identity in distribution
$$\zeta_{M(t)+1}\stackrel{d}{=}t+ (M(t)+1)\eta,$$
where $\eta$ is a unit exponential random variable, independent of $M(t)$. 
Thus
\begin{equation}\label{Ezeta}
\E \,\zeta_{M(t)+1}=t+\E\,{M(t)}+1,
\end{equation}
which together with   (\ref{Wald})   
give\footnote{For the general birth process as in the previous footnote, 
assuming $M(0)=0$ the identity
is
$$t=\E\left(\sum_{n=0}^{M(t)-1} \frac{1}{\beta_n}\right)$$}
\begin{equation}\label{fla}
\E{M^2(t)} =2t - \E {M(t)}.
\end{equation}

Alternatively, (\ref{fla}) can be derived from the $k=2$ instance of the formula
$$(\E M^k(t))'=\E \left( \frac{(M(t)+1)^k-M^k(t)}{M(t)+1}\right)$$
generalising (\ref{eqnnu}). Expanding the right-hand side, it is seen that  all moments ${\mathbb E}M^k(t)$ can be expressed, recursively, via the first moment  $\nu(s),~s\leq t$. 

Formula (\ref{fla})
allows us to express  the variance through  the mean as
\begin{equation}\label{varN}
\sigma^2(t) = 2t-\nu^2(t)-\nu(t).
\end{equation} 
Plugging the lower bound (\ref{bounds}) in (\ref{varN}) yields the bound
$\sigma^2(t) < \sqrt{2t+1}-1$, which for large $t$ is too far from (yet to be justified) (\ref{asympSigma}).
But working other way round we substitute (\ref{asympSigma}) with indefinite smaller order remainder in (\ref{varN}), 
and work out the quadratic equation 
 to extract the value of the sought  limit constant:
\begin{equation}\label{limconst}
\sqrt{2t}-\nu(t)\to \frac{2}{3}\,,~~~{\rm as}~t\to\infty.
\end{equation}
This result contrasts expansions (\ref{LIS}) and (\ref{LOIS}) but 
brings to mind some analogy with the expansion of the classic renewal function in the setting with uniform interarrival times \cite{Feller} (Ch. 11, Example 8).
Numerical 
calculations with (\ref{exact}) suggest that the limit is approached monotonically from below.

In a private communication, Alex Marynych informed  us that he could arrive at the asymptotics   $\sqrt{2t}-2/3$ using 
Theorem 1 from \cite{Minami} to approximate the tail of the Borel distribution in (\ref{Borel}).
Continuing the analogy with \cite{BD, GSBernoulli}, one can conjecture that the next term of the asymptotic expansion of the mean is of the order of $t^{-1/2}$, see also
\cite{Flaj} for a similar situation.

\section{The smallest first policy for fixed sample size }

We turn to the smallest first policy in the bin packing problem with unit capacity and fixed sample size $n$.
Let $0<u_{n1}<\dots<u_{nn}<1$ be the uniform $[0,1]$ order statistics, and let 
$$K_n=\max\{k: u_{n1}+\dots+u_{nk}\leq 1.$$ 
be the smallest first count, 
$\varkappa_n:={\mathbb E}K_n$. 
An explicit formula for $\varkappa_n$  exists  \cite{Coffman} (Theorem 7, with $c=1$), 
but is not  particularly user-friendly, as involving
an alternating double.
It is well known  that $\varkappa_n\sim \sqrt{2n}$   and that $\sqrt{2n}$ is also an upper bound  \cite{BR, Coffman, Gnedin, SS,  SteeleBR}

We assert now a much more precise result:
\begin{equation}\label{depois}
\varkappa_n=\sqrt{2n}-\frac{2}{3}+o(1), ~~~ n\to\infty.
\end{equation}
To show this, we first resort to the setting of the Poisson process on $[0,1]$ with rate $t$. The exit count $N(t)$ translates as the maximal number of Poisson points
whose total is at most $1$. Thus we have the poissonisation relation
$$\nu(t)=\sum_{n=1}^\infty  \varkappa_n   e^{-t}\frac{t^n}{n!}.$$

To depoissonise, we check conditions of Theorem 1 from \cite{Sz}.
The function $\nu(t)$ given by (\ref{exact}) for complex argument $t\in {\mathbb C}$ is an entire function, as the series converges everywhere. 
Some analytic work with the aid of the Stirling formula 
shows that $|\nu(t)|< c_1 |t|^{1/2}$ in the sector $|\arg t|<\pi/4$. 
Outside  the sector,  we have an estimate 
$$|e^t \nu(t)|<c_2 |t|^{1/2}\exp({|t|/\sqrt{2}}).$$
This follows by observing that 
 the maximum of    $| e^t (1-e^{-t/j})|$ for given $|t|$ is achieved at the boundary $|\!\arg t|=\pi/4$,
and then by approximating the sum (\ref{exact}) by an integral.
The cited theorem gives the possonisation error
$$|\nu(t)-\varkappa_{\lfloor t\rfloor}|= O(t^{-1/2}), ~~~~t\to\infty,$$
hence (\ref{limconst}) implies (\ref{depois}).

More generally, suppose the bin has capacity $C>0$. Consider the smallest first policy applied to the Poisson process on $[0,1]$ with rate $t$
and, in parallel, to $n$ items sampled from the  uniform $[0,1]$ distribution. Extending our notation from the case $C=1$, let  $N_C(t)$ and $K_{C,n}$ be the counts of items packed, and let $\nu_C(t), \varkappa_{C,n}$ be their means, respectively.

Generalising the $C=1$ result, we argue that
\begin{equation}\label{genC}
\nu_C(t)=\sqrt{2Ct}-\tfrac{2}{3}+o(1) {\rm ~~and~~} |\nu_C(t)-\varkappa_{C,\lfloor t \rfloor}|=O(t^{-1/2}).
\end{equation}

For $C\leq 1$, this is straightforward, as 
$\nu_C(t)=\nu(Ct)$ and we
readily  conclude (\ref{genC}) 
 from the $C=1$ case.

For $C>1$ we need to be more careful  since the maximum size of item is constrained by 1 and not by $C$.
Assessing the mean 
 in terms of the unit Poisson process on $[0, Ct]$ we have
\begin{eqnarray*}
\nu_C(t)=\sum_{n=1}^\infty \prob(N_C(t)\geq n)=
\sum_{n=1}^\infty \left\{
\prob(T(\pi_n)\leq Ct)-\prob(T(\pi_n)\leq Ct, \pi_n>t)\right\}=\\
\nu(Ct)-\sum_{n=1}^\infty \prob(T(\pi_{n-1})\leq (C-1)t, \pi_n>t).
\end{eqnarray*}
Recalling (\ref{TpiS}), we get the identity $(T(\pi_{n-1}), \pi_n)\stackrel{d}{=} (S(\pi_{n}), \pi_n)$, hence the last sum becomes
$$\sum_{n=1}^\infty \prob(S(\pi_{n})\leq (C-1)t, \pi_n>t)\leq\sum_{n=1}^{\lfloor t/2\rfloor}\prob(\pi_n>t)     +  \sum_{n={\lfloor t/2\rfloor}+1}^\infty \prob(S(\pi_n)\leq (C-1)t)=:\Sigma_1+\Sigma_2.$$
As $t\to\infty$, a large deviation  estimate for the Poisson process shows that $\Sigma_1$ approaches $0$ exponentially fast, and the same is shown for $\Sigma_2$ using 
$S(\pi_n)=\zeta_n$ and
the large deviation estimate (\ref{LaDe}) in the next section.

\section{A large deviation bound}

To justify  asymptotics of  the variance (\ref{asympSigma}) it remains to verify that the family
$$\frac{M(t)-\sqrt{2t}}{(\tfrac{1}{3}\sqrt{2t})^{1/2}}, ~~~t>0,$$
 is uniformly integrable.

We consider first
$${\zeta_n-a_n}=\sum_{j=1}^n j(\eta_j-1).$$  
For $\eta$ with the unit exponential distribution, the central moments are estimated as 
$$\E (\eta-1)^m=m!\sum_{j=1}^m\frac{(-1)^j}{j!}< \frac{m!}{e}+1.$$
Using this it is easy to check that 
$$
\E[k(\eta-1)]^m\leq \frac{k^2}{2} n^{m-2} m!\,,~~~~~1\leq k\leq n,~m\geq 2,
$$
which verifies the  condition for  large deviation bounds from  \cite{Petrov} (Chapter 3, Theorem 17). Hence
we obtain 
\begin{equation}\label{LaDe}
\sup_{{z\geq 0}}\prob(|\zeta_n-a_n|> b_n z)< 3\,e^{-z/4},
\end{equation}
where both constants  are not sharp.
Inverting the latter we arrive at similar bound
$$
\prob\left(|M(t)-\sqrt{2t}|>z \left(\tfrac{1}{3}\sqrt{2t}\right)^{1/2}\right)< c \,e^{-z^2/4},
$$
which implies the desired uniform integrability.

\end{document}